\documentclass{amsart}

\usepackage[all]{xypic}
\usepackage[centertags]{amsmath}
\usepackage{amsfonts}
\usepackage{amscd}
\usepackage{amssymb}
\usepackage{amsthm}
\usepackage{newlfont}
\usepackage{amsxtra}
 \newtheorem{thm}{Theorem}[section]
 \newtheorem{cor}[thm]{Corollary}
 
 \newtheorem{lem}[thm]{Lemma}
 \newtheorem{prop}[thm]{Proposition}
 \theoremstyle{definition}
 \newtheorem{defn}[thm]{Definition}
 \theoremstyle{remark}
 \newtheorem{rem}[thm]{Remark}
 \theoremstyle{remark}
 \newtheorem{example}[thm]{Example}
 \theoremstyle{definition}
 \newtheorem{notn}[thm]{Notation}
 \numberwithin{equation}{section}

 \newcommand{\an}{\mathrm{an}}

 \newcommand{\Ver}{\mathrm{Ver}}
 
 \newcommand{\Spec}{\mathrm{Spec}}

 \newcommand{\Aut}{\mathrm{Aut}}
 
 \newcommand{\Pic}{\mathrm{Pic}}
 
 \newcommand{\ord}{\mathrm{ord}}

 \newcommand{\GL}{\mathrm{GL}}
 \newcommand{\PGL}{\mathrm{PGL}}
 \newcommand{\SL}{\mathrm{SL}}

 \newcommand{\tor}{\mathrm{tor}}

 \newcommand{\Ed}{\mathrm{Ed}}

 \newcommand{\Tr}{\mathrm{Tr}}

 \renewcommand{\mod}{\mathrm{mod}}
 \newcommand{\Nr}{\mathrm{Nr}}

 \newcommand{\fp}{\mathfrak p}
 \newcommand{\fr}{\mathfrak r}

 \newcommand{\cO}{\mathcal{O}}

 \newcommand{\cA}{\mathcal{A}}
 
 \newcommand{\cG}{\mathcal{G}}
 \newcommand{\cP}{\mathcal{P}}

 \renewcommand{\cD}{\mathcal{D}}
 
 \newcommand{\cX}{\mathcal{X}}

 \newcommand{\cI}{\mathcal{I}}
 
 \renewcommand{\cH}{\mathcal{H}}
 \newcommand{\cT}{\mathcal{T}}

 \newcommand{\R}{\mathbb{R}}

 \newcommand{\F}{\mathbb{F}}
 \newcommand{\M}{\mathbb{M}}
 \newcommand{\Q}{\mathbb{Q}}
 \newcommand{\Z}{\mathbb{Z}}
 \newcommand{\A}{\mathbb{A}}

 \renewcommand{\P}{\mathbb{P}}

 \newcommand{\bs}{\setminus}

 \newcommand{\Fi}{F_\infty}

 \newcommand{\G}{\Gamma}
 
 \newcommand{\La}{\Lambda}
 \newcommand{\la}{\lambda}

\begin{document}

\title{Trees, quaternion algebras and modular curves}

\author{Mihran Papikian}

\address{Department of Mathematics, Pennsylvania State University, University Park, PA 16802}

\email{papikian@math.psu.edu}

\thanks{The author was supported in part by NSF grant DMS-0801208 and Humboldt Research Fellowship.}



\begin{abstract}
We study the action on the Bruhat-Tits tree of unit groups of
maximal orders in certain quaternion algebras over $\F_q(T)$ and
discuss applications to arithmetic geometry and group theory.
\end{abstract}


\maketitle


\section{Introduction}\label{SecIntr}

Let $B$ be an indefinite division quaternion algebra over $\Q$. Let
$\cO\subset B$ be a maximal order. Let $d$ be the discriminant of
$B$; recall that $d>1$ is the product of primes where $B$ ramifies.
Let $\G^d=\{\gamma\in \cO\ |\ \Nr(\gamma)=1\}$, where $\Nr$ is the
reduced norm of $B$. Upon fixing an identification of $B\otimes_\Q
\R$ with $\M_2(\R)$, we can view the group $\G^d$ as a discrete
subgroup of $\SL_2(\R)$. It is well-known that $\G^d$ acts
discontinuously on the upper half-plane $\cH$ and $X^d:=\G^d\bs \cH$
is compact; cf. \cite[Ch. 5]{Miyake}. Shimura observed that $X^d$
parametrizes abelian surfaces with multiplication by $\cO$. From
this he was able to show that the algebraic curve $X^d$ admits a
canonical model over $\Q$. It is hardly surprising that the
properties of the algebraic curve $X^d$ are reflected in the
properties of the group $\G^d$. For example, one can compute the
genus of $X^d$ by computing the hyperbolic volume of a fundamental
domain for the action of $\G^d$ on $\cH$ and the number of elliptic
points; cf. \cite[Ch. 4]{Vigneras}. The transfer of information also
goes in the opposite direction, e.g., knowing the genus of $X^d$ and
the relative position of elliptic points on this curve, one deduces
a presentation for $\G^d$. An interesting special case is when
$\G^d$ is generated by torsion elements; this happens only for
$d=6,10,22$, and one knows a presentation of $\G^d$ in these cases,
cf. \cite{AB}.

\vspace{0.1in}

The purpose of the present paper is to study the function field
analogues of the groups $\G^d$. First, we introduce some notation.

Let $C:=\P^1_{\F_q}$ be the projective line over the finite field
$\F_q$. Denote by $F=\F_q(T)$ the field of rational functions on
$C$. The set of closed points on $C$ (equivalently, places of $F$)
is denoted by $|C|$. For each $x\in |C|$, we denote by $\cO_x$ and
$F_x$ the completions of $\cO_{C,x}$ and $F$ at $x$, respectively.
The residue field of $\cO_x$ is denoted by $\F_x$, the cardinality
of $\F_x$ is denoted by $q_x$, and $\deg(x):=[\F_x:\F_q]$. Let
$A:=\F_q[T]$ be the ring of polynomials in $T$ with $\F_q$
coefficients; this is the subring of $F$ consisting of functions
which are regular away from $\infty:=1/T$.

Let $D$ be a quaternion division algebra over $F$. Let $R$ be the
set of places where $D$ ramifies ($R$ is finite and has even
cardinality, and conversely, for any choice of a finite set
$R\subset |C|$ of even cardinality there is a unique, up to an
isomorphism, quaternion algebra ramified exactly at the places in
$R$; see \cite[p. 74]{Vigneras}). Denote by $D^\times$ the
multiplicative group of $D$. Assume $D$ is split at $\infty$, i.e.,
$D\otimes_F\Fi\cong \M_2(\Fi)$. Fix a maximal $A$-order $\La$ in
$D$. Since $D$ is split at $\infty$, it satisfies the so-called
\textit{Eichler condition} relative to $A$. Since $\Pic(A)=1$, this
implies that, up to conjugation, $\La$ is the unique maximal
$A$-order in $D$, i.e., any other maximal $A$-order in $D$ is of the
form $\alpha \La\alpha^{-1}$ for some $\alpha\in D^\times(F)$, cf.
\cite[Cor. III. 5.7]{Vigneras}. We are interested in the group of
units of $\La$:
$$
\G:=\La^\times=\{\la\in \La\ |\ \Nr(\la)\in \F_q^\times\}.
$$

Via an isomorphism $D^\times(\Fi)\cong \GL_2(\Fi)$, the group $\G$
can be considered as a discrete subgroup of $\GL_2(\Fi)$. There are
two analogues of the Poincar\'e upper half-plane in this setting.
One is the \textit{Bruhat-Tits tree} $\cT$ of $\PGL_2(\Fi)$, the
other is \textit{Drinfeld upper half-plane}
$\Omega:=\P^{1,\an}_{\Fi}-\P^{1,\an}_{\Fi}(\Fi)$, where
$\P^{1,\an}_{\Fi}$ be the rigid-analytic space associated to the
projective line over $\Fi$. These two versions of the upper
half-plane are related to each other: $\Omega$ has a natural
structure of a smooth geometrically connected rigid-analytic space
and $\cT$ is the dual graph of an analytic reduction of $\Omega$,
cf. \cite{vdPut}. As a subgroup of $\GL_2(\Fi)$, $\G$ acts naturally
on both $\cT$ and $\Omega$ (these actions are compatible with
respect to the reduction map). The quotient $\G\bs \Omega$ is a
one-dimensional, connected, smooth analytic space over $\Fi$, which
in fact is the rigid-analytic space associated to a smooth,
projective curve $X^R$ over $\Fi$, cf. \cite[Thm. 3.3]{vdPut}. The
quotient $\G\bs \cT$ is a finite graph (Proposition \ref{propFG}).

\vspace{0.1in}

The properties of $\G$, $\G\bs\cT$ and $X^R$ are intimately related.
We explore these relationships to obtain a description of the graph
$\G\bs \cT$, which includes a formula for the number of vertices of
a given degree and a formula for the first Betti number (Theorem
\ref{PropTV}). Next, the description of $\G\bs \cT$ translates into
a statement about the structure of the group $\G$ (Theorem
\ref{thmGTI}). In particular, we determine the cases when $\G$ is
generated by torsion elements and find a presentation for $\G$ in
those cases (Theorem \ref{thm-tree}); in $\S$\ref{SecExplU}, we find
explicitly the torsion units which generate $\G$ in terms of a basis
of $D$. Finally, we have an application to the arithmetic of $X^R$.
From the structure of $\G\bs \cT$ and a geometric version of
Hensel's lemma we deduce that $X^R(\Fi)=\emptyset$ if and only if
$R$ contains a place of even degree (Theorem \ref{thmFipoints}).
This last theorem is the function field analogue of a well-known
result of Shimura which says that $X^d(\R)=\emptyset$, i.e., $X^d$
does not have points which are rational over $\R$ \cite{Shimura}.

The curve $X^R$ parametrizes $\cD$-elliptic sheaves over $F$ with
pole $\infty$ \cite{BS}, where $\cD$ is a maximal $\cO_C$-order in
$D$ (the notion $\cD$-elliptic sheaf generalizes the notion of
Drinfeld module). This can be used to show that $X^R$ has a model
over $F$. As a consequence of Theorem \ref{thmFipoints}, if $R$
contains a place of even degree, then $X^R$ has no $L$-rational
points for any extension $L/F$ which embeds into $\Fi$. In
particular, $X^R(F)=\emptyset$.


\section{Preliminaries}

\subsection{Graphs}\label{SecG} We recall the terminology related to graphs, as presented in
\cite{SerreT}.

\begin{defn}
An (oriented) \textit{graph} $\cG$ consists of a non-empty set
$X=\Ver(\cG)$, a set $Y=\Ed(\cG)$ and two maps
$$
Y\to X\times X, \quad y\mapsto (o(y), t(y))
$$
and
$$
Y\to Y, \quad y\mapsto \bar{y}
$$
which satisfy the following condition: for each $y\in Y$ we have
$\bar{\bar{y}}=y$, $\bar{y}\neq y$ and $o(y)=t(\bar{y})$.
\end{defn}

An element $v\in X$ is called a \textit{vertex} of $\cG$; an element
$y\in Y$ is called an (\textit{oriented}) \textit{edge}, and
$\bar{y}$ is called the \textit{inverse} of $y$. The vertices $o(y)$
and $t(y)$ are the \textit{origin} and the \textit{terminus} of $y$,
respectively. These two vertices are called the \textit{extremities}
of $y$. Note that it is allowed for distinct edges $y\neq z$ to have
$o(y)=o(z)$ and $t(y)=t(z)$:
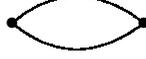
\begin{center}
\begin{picture}(35,20)
\qbezier(5,10)(17,20)(30,10)\qbezier(5,10)(17,0)(30,10)
\put(5,10){\circle*{2}}\put(30,10){\circle*{2}}
\end{picture}
\end{center}
and it is also allowed to have $y\in Y$ with $o(y)=t(y)$, in which
case $y$ is called a \textit{loop}:
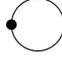
\begin{center}
\begin{picture}(20,20)
\put(10,10){\circle{10}}\put(5,10){\circle*{2}}
\end{picture}
\end{center}
We say that two vertices are \textit{adjacent} if they are the
extremities of some edge. We will assume that for any $v\in X$ the
number of edges $y\in Y$ with $o(y)=v$ is finite; this number is the
\textit{degree} of $v$. A vertex $v\in X$ is called
\textit{terminal} if it has degree $1$. A \textit{circuit} in $\cG$
is a collection of edges $y_1,y_2,\dots, y_m\in Y$ such that
$y_i\neq \bar{y}_{i+1}$, $t(y_i)=o(y_{i+1})$ for $1\leq i\leq m-1$,
and $t(y_m)=o(y_1)$. A graph $\cG$ is \textit{connected} if for any
two distinct vertices $v, w\in X$ there is a collection of edges
$y_1,y_2,\dots, y_m\in Y$ such that $t(y_i)=o(y_{i+1})$ for $1\leq
i\leq m-1$, $v=o(y_1)$ and $w=t(y_m)$.

A connected graph without circuits is called a \textit{tree}. A
\textit{geodesic} between two vertices $v, w$ in a tree $\cT$ is a
collection of edges $y_1,y_2,\dots, y_m\in \Ed(\cT)$ such that
$y_i\neq \bar{y}_{i+1}$, $t(y_i)=o(y_{i+1})$ for $1\leq i\leq m-1$,
and $v=o(y_1)$, $w=t(y_m)$. Since $\cT$ is connected and has no
circuits, a geodesic between $v, w\in \Ver(\cT)$ always exists and
is unique. The distance $d(v, w)$ between $v$ and $w$ is the number
of edges in the geodesic joining $v$ and $w$. In particular, $d(v,
w)=1$ if and only if $v$ and $w$ are adjacent.

A graph is \textit{finite} if it has finitely many vertices and
edges. A finite graph $\cG$ can be interpreted as a $1$-dimensional
$\Delta$-complex \cite[p. 102]{Hatcher}. The \textit{first Betti
number} $h_1(\cG)$ of $\cG$ is the dimension of the homology group
$\dim_\Q H_1(\cG, \Q)$. (One can show that $\cG$ is homotopic to a
bouquet of circles; $h_1(\cG)$ is equal to the number of those
circles.)

An \textit{automorphism} of $\cG$ is a pair $\phi=(\phi_1, \phi_2)$
of bijections $\phi_1: X\to X$ and $\phi_2: Y\to Y$ such that
$\phi_1(o(y))=o(\phi_2(y))$ and
$\overline{\phi_2(y)}=\phi_2(\bar{y})$. Let $\G$ be a group acting
on a graph $\cG$ (i.e., $\G$ acts via automorphisms). We say that
$v, w\in X$ are $\G$-\textit{equivalent} if there is $\gamma \in \G$
such that $\gamma v=w$; similarly, $y, z\in Y$ are
$\G$-\textit{equivalent} if there is $\gamma\in \G$ such that
$\gamma y=z$. For $v\in X$, denote
$$
\G_v=\{\gamma\in G\ |\ \gamma v=v\}
$$
the stabilizer of $v$. Similarly, let $\G_{y}=\G_{\bar{y}}$ be the
stabilizer of $y\in Y$. $\G$ acts with \textit{inversion} if there
is $\gamma\in \G$ and $y\in Y$ such that $\gamma y=\bar{y}$. If $\G$
acts without inversion, then we have a natural quotient graph
$\G\bs\cG$ such that $\Ver(\G\bs \cG)=\G\bs \Ver(\cG)$ and
$\Ed(\G\bs \cG)=\G\bs \Ed(\cG)$.

\begin{defn}(\cite[p.70]{SerreT}) Let $K$ be a field complete with respect to a
discrete valuation $\ord_K$. Let $\cO:=\{x\in K\ |\ \ord_K(x)\geq
0\}$ be the ring of integers of $\cO$. Fix a uniformizer $\pi$,
i.e., an element with $\ord_K(\pi)=1$. Let $k:=\cO/\pi\cO$ be the
residue field. We assume $k\cong \F_q$ is finite. Let $V$ be a
two-dimensional vector space over $K$. A \textit{lattice} of $V$ is
any finitely generated $\cO$-submodule of $V$ which generates the
$K$-vector space $V$; such a module is free of rank $2$. If $\La$ is
a lattice and $x\in K^\times$, then $x\La$ is also a lattice. We
call $\La$ and $x\La$ equivalent lattices. The action of $K^\times$
on the set of lattices in $V$ subdivides this set into disjoint
equivalence classes. We denote the class of $\La$ by $[\La]$.

Let $\cT$ be the graph whose vertices $\Ver(\cT)=\{[\La]\}$ are the
equivalence classes of lattices in $V$, and two vertices $[\La]$ and
$[\La']$ are adjacent if we can choose representatives $L\in [\La]$
and $L'\in [\La']$ such that $L'\subset L$ and $L/L'\cong k$. One
shows that $\cT$ is an infinite tree in which every vertex has
degree $(q+1)$. This is the \textit{Bruhat-Tits tree} of
$\PGL_2(K)$.
\end{defn}

Let $\GL(V)$ denote the group of $K$-automorphisms of $V$; it is
isomorphic to $\GL_2(K)$. One easily verifies that $\GL(V)$ acts on
the Bruhat-Tits tree $\cT$ and preserves the distance between any
two vertices (although $\GL(V)$ acts with inversion). Let
$$
\GL(V)^0:=\ker\left(\ord_K\circ \det: \GL(V)\to \Z\right).
$$
\begin{lem}\label{lem1.6}
Let $g\in \GL(V)^0$ and $v\in \Ver(\cT)$. Then $d(v, gv)=2n$, $n\geq
0$.
\end{lem}
\begin{proof}
See \cite[p. 75]{SerreT}.
\end{proof}

Given two linearly independent vectors $f_1$ and $f_2$ in $V$, we
denote by $[f_1,f_2]$ the similarity class of the lattice $\cO
f_1\oplus \cO f_2$. Fix the standard basis $e_1=(1,0)$, $e_2=(0,1)$
of $K^2$. Note that the vertices adjacent to $[e_1, e_2]$ are
$$
[\pi^{-1}e_1, e_2], \quad [\pi e_1, e_2+c e_1], c\in k.
$$
If we express every vector in $K^2$ in terms of $e_1$ and $e_2$,
then the action of $g\in \GL_2(K)$ on $\cT$ is explicitly given by
$$
g[ae_1+be_2, ce_1+de_2]=[ae_1g+be_2g, ce_1g+de_2g],
$$
where $g$ acts on row vectors $e_i$ in the usual manner.

\subsection{Quaternion algebras}\label{Sec3}

Given a quaternion algebra $D$ over $F$, we denote by $\alpha\mapsto
\alpha'$ the canonical involution of $D$ \cite[p. 1]{Vigneras}; thus
$\alpha''=\alpha$ and $(\alpha\beta)'=\beta'\alpha'$. The
\textit{reduced trace} of $\alpha$ is $\Tr(\alpha)=\alpha+\alpha'$;
the \textit{reduced norm} of $\alpha$ is
$\Nr(\alpha)=\alpha\alpha'$; the \textit{reduced characteristic
polynomial} of $\alpha$ is
$$f(x)=(x-\alpha)(x-\alpha')=x^2-\Tr(\alpha)x+\Nr(\alpha).$$

\begin{notn} For $a, b\in F^\times$, let $H(a,b)$ be the $F$-algebra with basis $1,i,j,
ij$ (as an $F$-vector space), where $i,j\in H(a,b)$ satisfy
\begin{itemize}
\item If $q$ is odd,
$$
i^2=a,\quad j^2=b,\quad ij=-ji;
$$
\item If $q$ is even,
$$
i^2+i=a,\quad j^2=b,\quad ij=j(i+1).
$$
\end{itemize}
\end{notn}

\begin{prop}
$H(a,b)$ is a quaternion algebra. Moreover, any quaternion algebra
$D$ is isomorphic to $H(a,b)$ for some $a,b\in F^\times$ $($although
$a$ and $b$ are not uniquely determined by $D$$)$.
\end{prop}
\begin{proof}
See pp. 1-5 in \cite{Vigneras}.
\end{proof}

For $\alpha\in H(a,b)$, there is a unique expression
$\alpha=x+yi+zj+wij$. The \textit{quadratic form} corresponding to
the reduced norm on $H(a,b)$ is $Q(x,y,z,w)=\alpha\cdot \alpha'$. In
terms of $a$ and $b$, $Q$ is explicitly given by
\begin{itemize}
\item If $q$ is odd,
$$
Q(x,y,z,w)=x^2-ay^2-bz^2+abw^2;
$$
\item If $q$ is even,
$$
Q(x,y,z,w)=x^2+xy+ay^2+b(z^2+zw+aw^2).
$$
\end{itemize}

One has the following useful criterion for determining whether
$H(a,b)$ is split or ramified at a given place.

\begin{lem}\label{lem-split}
$H(a,b)$ is split at $v\in |C|$ if and only if $Q(x,y,z,w)=0$ has a
non-trivial solution over $F_v$, i.e., there exist $(c_1, c_2,
c_3,c_4)\neq (0,0,0,0)$ such that $Q(c_1, c_2, c_3,c_4)=0$.
\end{lem}
\begin{proof}
Let $0\neq \alpha\in H(a,b)_v:=H(a,b)\otimes_F F_v$. If $Q$ does not
have non-trivial solutions over $F_v$, then $\alpha\cdot \alpha'\in
F_v^\times$ and $\alpha'/(\alpha\cdot \alpha')$ is the inverse of
$\alpha$. On the other hand, if $Q$ does have non-trivial solutions
then there is $\alpha\neq 0$ such that $\alpha\cdot \alpha'=0$. In
this last case $H(a,b)_v$ obviously has zero divisors, so cannot be
a division algebra.
\end{proof}

\begin{cor}\label{cor-split2}
If $a,b\in \F_q^\times$, then $H(a,b)$ is split, i.e., is isomorphic
to $\M_2(F)$.
\end{cor}
\begin{proof}
Using Lemma \ref{lem-split} and Hensel's lemma, $H(a,b)$ is split at
$v$ if and only if $Q$ has a non-trivial solution over $\F_v$. An
easy consequence of Warning's Theorem \cite[Thm. 6.5]{LN} is that a
quadratic form in $n$ variables over a finite field has a
non-trivial solution if $n\geq 3$. Since $Q$ over $\F_v$ has $4$
variables, it has a non-trivial solution. Hence $H(a,b)$ is split at
every $v\in |C|$; this implies $H(a,b)\cong \M_2(F)$ \cite[p.
74]{Vigneras}.
\end{proof}

\begin{prop}\label{prop-emb}
Let $L$ be a finite field extension of $F$. There is an
$F$-isomorphism of $L$ onto an $F$-subalgebra of $D$ if and only if
$[L:F]$ divides $2$ and no place in $R$ splits in $L$. Moreover, any
two such $F$-isomorphisms are conjugate in $D$.
\end{prop}
\begin{proof}
See \cite[(32.15)]{Reiner}.
\end{proof}


\section{The group of units}\label{SecGU} We return to the notation
and assumptions in the introduction. Thus, $D$ is a quaternion
division algebra which is split at $\infty$, $\La$ is a maximal
$A$-order in $D$, and
$$
\G:=\La^\times=\{\la\in \La\ |\ \Nr(\la)\in \F_q^\times\}.
$$

From now on we denote $K:=\Fi$, $\cO:=\cO_\infty$,
$k:=\F_\infty\cong\F_q$. The group $\G$ can be considered as a
discrete subgroup of $\GL_2(K)$ via an embedding
$$
\iota:\G\hookrightarrow D^\times(F)\hookrightarrow D^\times(K)\cong
\GL_2(K).
$$
Note that for $\gamma\in \G$, $\det(\iota(\gamma))=\Nr(\gamma)$, so
$\iota(\Gamma)\subset \GL_2(K)^0$. We fix some embedding $\iota$,
and omit it from notation. Being a subgroup of $\GL_2(K)^0$, $\G$
naturally acts on the Bruhat-Tits tree $\cT$ of $\PGL_2(K)$, and
moreover, by Lemma \ref{lem1.6}, $\G$ acts without inversion.

\begin{prop}\label{propFG}
The quotient graph $\G\bs \cT$ is finite.
\end{prop}
\begin{proof}
It is enough to show that $\G\bs \cT$ has finitely many vertices.
The group $\GL_2(K)$ acts transitively on the set of lattices in
$K^2$. The stabilizer of $[e_1, e_2]$ is $Z(K)\cdot \GL_2(\cO)$,
where $Z$ denotes the center of $\GL(2)$. This yields a natural
bijection $\Ver(\cT)\cong \GL_2(K)/Z(K)\cdot \GL_2(\cO)$, and also
$$
\Ver(\G\bs \cT)\cong \G\bs \GL_2(K)/Z(K)\cdot \GL_2(\cO).
$$
We will show that the above double coset space is finite. Denote by
$\A$ the adele ring of $F$. Consider the group $D^\times(F)$
embedded diagonally into $D^\times(\A)$. Since $D$ is a division
algebra, $D^\times(F)\bs D^\times(\A)/Z(K)$ is compact (cf.
\cite[Ch. III.1]{Vigneras}). Denote $\cD_f^\times:=\prod_{x\in
|C|-\infty}D^\times(\cO_x)$. Since $A$ is a principal ideal domain,
the strong approximation theorem for $D^\times$ yields (cf.
\cite[Ch. III.4]{Vigneras})
$$
D^\times(\A)\cong D^\times(F)\cdot \GL_2(K)\cdot \cD_f^\times.
$$
Note that $\G=D^\times(F)\cap \cD_f^\times$. Thus, $\G\bs
\GL_2(K)/Z(K)$ is compact since it is the image of $D^\times(F)\bs
D^\times(\A)/Z(K)$ under the natural quotient map $D^\times(\A)\to
D^\times(\A)/\cD^\times_f$. Finally, since $\GL_2(\cO)$ is open in
$\GL_2(K)$, $\G\bs \GL_2(K)/Z(K)\cdot \GL_2(\cO)$ is finite.
\end{proof}

\begin{prop}\label{Prop2.1}
Let $v\in \Ver(\cT)$ and $y\in \Ed(\cT)$. Then $\G_v\cong
\F_q^\times$ or $\G_v\cong \F_{q^2}^\times$, and $\G_y\cong
\F_q^\times$.
\end{prop}
\begin{proof}
By choosing an appropriate basis of $K^2$, we can assume $v=[e_1,
e_2]$. The stabilizer of $v$ in $\GL_2(K)^0$ is $\GL_2(\cO)$. Since
this last group is compact in $\GL_2(K)$, whereas $\G$ is discrete,
$\G_v=\G\cap \GL_2(\cO)$ is finite. In particular, if $\gamma\in
\G_v$, then $\gamma^n=1$ for some $n\geq 1$. We claim that the order
$n$ of $\gamma$ is coprime to the characteristic $p$ of $F$. Indeed,
if $p|n$ then $(\gamma^{n/p}-1)\in D$ is non-zero but
$(\gamma^{n/p}-1)^p=0$. This is not possible since $D$ is a division
algebra. Consider the subfield $F(\gamma)$ of $D$ generated by
$\gamma$ over $F$. By Proposition \ref{prop-emb}, $[F(\gamma):F]=1$
or $2$. Since $\gamma\in D$ is algebraic over $\F_q$, we conclude
that $[\F_q(\gamma):\F_q]=1$ or $2$.

It is obvious that $\F_q^\times\subset \G_v$. Assume there is
$\gamma\in \G_v$ which is not in $\F_q^\times$. From the previous
paragraph, $\gamma$ generates $\F_{q^2}$ over $\F_q$. Considering
$\gamma$ as an element of $\GL_2(\cO)$, we clearly have $a+b
\gamma\in \M_2(\cO)$ for $a, b\in \F_q$ (embedded diagonally into
$\GL_2(K)$). But if $a$ and $b$ are not both zero, then
$a+b\gamma\in \La$ is invertible, hence belongs to $\G$ and
$\GL_2(\cO)$. We conclude that
$\F_q(\gamma)^\times\cong\F_{q^2}^\times\subset \G_v$, and moreover,
every element of $\G_v$ is of order dividing $q^2-1$. Suppose there
is $\delta\in \G_v$ which is not in $\F_q(\gamma)^\times$. Since
$\delta$ is algebraic over $\F_q$, $\delta$ and $\gamma$ do not
commute in $D$ (otherwise $F(\gamma, \delta)$ is a subfield of $D$
of degree $>2$ over $F$). Then $\G_v/\F_q^\times$ is a finite
subgroup of $\PGL_2(K)$ whose elements have orders dividing $q+1$
and which contains two non-commuting elements of order $(q+1)$. This
contradicts Dickson's classification of finite subgroups of
$\PGL_2(K)$ \cite[II.8.27]{Huppert}.

Now consider $\G_y$. Clearly $\F_q^\times\subset \G_y$. Let $v$ and
$w$ be the extremities of $y$. Note that there are natural
inclusions $\G_y\subset \G_v$, $\G_y\subset \G_w$ and $\G_y=\G_v\cap
\G_w$. If $\G_y$ is strictly larger than $\F_q^\times$, then from
the discussion about the stabilizers of vertices, we have
$\G_v=\G_w\cong \F_{q^2}^\times$ (an equality of subgroups of $\G$).
Therefore, $\G_y\cong \F_{q^2}^\times$. On the other hand, the
stabilizer of $y$ in $\GL_2(K)^0$ is isomorphic to the Iwahori
subgroup $\cI$ of $\GL_2(\cO)$ consisting of matrices $\begin{pmatrix} a & b \\
c & d
\end{pmatrix}$ such that $c\equiv 0\ (\mod\ \pi)$. Since $\cI$ does
not contain a subgroup isomorphic to $\F_{q^2}^\times$, we get a
contradiction.
\end{proof}

\begin{cor}\label{cor2.2} Let $v\in \Ver(\cT)$ be such that
$\G_v\cong \F_{q^2}^\times$. Then $\G_v$ acts transitively on the
vertices adjacent to $v$.
\end{cor}
\begin{proof}
By Proposition \ref{Prop2.1}, a subgroup of $\F_{q^2}^\times$ which
stabilizes an edge with origin $v$ is $\F_q^\times$. Hence
$\G_v/\F_q^\times$ acts freely on the set of vertices adjacent to
$v$. Since this quotient group has $q+1$ elements, which is also the
number of vertices adjacent to $v$, it has to act transitively.
\end{proof}

We introduce a function which will simplify the notation in our
later discussions:
$$
\wp(R)=\left\{
         \begin{array}{ll}
           0, & \hbox{if some place in $R$ has even degree;} \\
           1, & \hbox{otherwise.}
         \end{array}
       \right.
$$
Let
$$
g(R)=1+\frac{1}{q^2-1}\prod_{x\in R}(q_x-1)-\frac{q}{q+1}\cdot 2^{\#
R-1}\cdot \wp(R).
$$

\begin{thm}\label{PropTV}\hfill
\begin{enumerate}
\item The graph $\G\bs \cT$ has no loops;
\item $h_1(\G\bs \cT)=g(R)$;
\item Every vertex of $\G\bs \cT$ is either terminal or has degree
$q+1$.
\item Let $V_1$ and $V_{q+1}$ be the number of terminal and
degree $q+1$ vertices of $\G\bs \cT$, respectively. Then
$$
V_1=2^{\# R -1}\wp(R)\quad \text{and}\quad
V_{q+1}=\frac{1}{q-1}(2g(R)-2+V_1).
$$
\end{enumerate}
\end{thm}
\begin{proof} The graph $\G\bs \cT$ has no loops since adjacent
vertices of $\cT$ are not $\G$-equivalent, as follows from Lemma
\ref{lem1.6}.

By Theorem \ref{thmKM}, $\G\bs \cT$ is the dual graph of
$\cX^R_{\bar{k}}$ (see $\S$\ref{SecMC} for notation). Using
\cite[Exr. IV.1.8]{Hartshorne}, one concludes that the arithmetic
genus of $\cX^R_{\bar{k}}$ is equal to $h_1(\G\bs \cT)$. Now by the
flatness of $\cX^R$, the arithmetic genus of $\cX^R_{\bar{k}}$ is
equal to the genus of $\cX^R_K\cong X^R$. Finally, the genus of
$X^R$ is equal to $g(R)$ by \cite[Thm. 5.4]{PapGenus}. Overall, we
get $h_1(\G\bs \cT)=g(R)$.

Let $v\in \Ver(\cT)$. By Proposition \ref{Prop2.1}, $\G_v\cong
\F_q^\times$ or $\F_{q^2}^\times$. In the second case, by Corollary
\ref{cor2.2}, the image of $v$ in $\G\bs \cT$ is a terminal vertex.
Now assume $\G_v\cong \F_q^\times$. We claim that the image of $v$
in $\G\bs \cT$ has degree $q+1$. Let $e, y\in \Ed(\cT)$ be two
distinct edges with origin $v$. It is enough to show that $e$ is not
$\G$-equivalent to $y$ or $\bar{y}$. On the one hand, $e$ cannot be
$\G$-equivalent to $y$ since $\G_v\cong \F_q^\times$ stabilizes
every edge with origin $v$. On the other hand, if $e$ is
$\G$-equivalent to $\bar{y}$ then $v$ is $\G$-equivalent to an
adjacent vertex, and that cannot happen.

Let $S$ be the set of terminal vertices of $\G\bs \cT$. Let $G$ be
the set of conjugacy classes of subgroups of $\G$ isomorphic to
$\F_{q^2}^\times$. We claim that there is a bijection $\varphi:S\to
G$ given by $\tilde{v}\mapsto \G_v$, where $v$ is a preimage of the
terminal vertex $\tilde{v}\in S$. The map is well-defined since if
$w$ is another preimage of $\tilde{v}$ then $v=\gamma w$ for some
$\gamma\in \G$, and so $\G_w=\gamma^{-1} \G_v\gamma$ is a conjugate
of $\G_v$. If $\varphi$ is not injective, then there are two
vertices $v,w\in\Ver(\cT)$ such that $\G_v\cong\F_{q^2}^\times$,
$\G_w=\gamma^{-1}\G_v\gamma$ for some $\gamma\in \G$, but $v$ and
$w$ are not in the same $\G$-orbit. Then $\G_{\gamma
w}=\gamma\G_w\gamma^{-1}=\G_v$, but $\gamma w\neq v$. The geodesic
connecting $v$ to $\gamma w$ is fixed by $\G_v$, so every edge on
this geodesic has stabilizer equal to $\G_v$. This contradicts
Proposition \ref{Prop2.1}. Finally, to see that $\varphi$ is
surjective it is enough to show that every torsion element in $\G$
fixes some vertex. Suppose $\gamma\in \G$ does not fix any vertices
in $\cT$. Since $\gamma$ acts without inversion, a result of Tits
\cite[Prop. 24, p. 63]{SerreT} implies that there is a straight path
$\cP$ in $\cT$ on which $\gamma$ induces a translation of amplitude
$m\geq 1$. If $\gamma$ is torsion, then for some $n\geq 2$
$\gamma^n=1$ induces a translation of amplitude $mn$ on $\cP$, which
is absurd.

Let $\cA:=\F_{q^2}[T]$ and $L:=\F_{q^2}F$ (note that $\cA$ is the
integral closure of $A$ in $L$). If $q$ is odd, let $\xi$ be a fixed
non-square in $\F_q$. If $q$ is even, let $\xi$ be a fixed element
of $\F_q$ such that $\Tr_{\F_q/\F_2}(\xi)=1$ (such $\xi$ always
exists, cf. \cite[Thm. 2.24]{LN}). Consider the polynomial
$f(x)=x^2-\xi$ if $q$ is odd, and $f(x)=x^2+x+\xi$ if $q$ is even.
Note that $f(x)$ is irreducible over $\F_q$; this is obvious for $q$
odd, and follows from \cite[Cor. 3.79]{LN} for $q$ even. Thus, a
solution of $f(x)=0$ generates $\F_{q^2}$ over $\F_q$. Denote by $P$
the set of $\G$-conjugacy classes of elements of $\La$ whose reduced
characteristic polynomial is $f(x)$. By a theorem of Eichler (Cor.
5.12, 5.13, 5.14 on pp. 94-96 of \cite{Vigneras}),
$$
\# P = h(\cA)\prod_{x\in R}\left(1-\left(\frac{L}{x}\right)\right),
$$
where $h(\cA)$ is the class number of $\cA$ and
$\left(\frac{L}{x}\right)$ is the \textit{Artin-Legendre symbol}:
$$
\left(\frac{L}{x}\right)=\left\{
                           \begin{array}{ll}
                             1, & \hbox{if $x$ splits in $L$;} \\
                             -1, & \hbox{if $x$ is inert in $L$;} \\
                             0, & \hbox{if $x$ ramifies in $L$.}
                           \end{array}
                         \right.
$$
Note that a place of even degree splits in $L$, and a place of odd
degree remains inert, and since $h(\cA)=1$, this implies $\# P =
2^{\# R}\wp(R)$.

If $\la\in \La$ is an element with reduced characteristic polynomial
$f(x)$, then it is clear that $\la\in \G$ and
$\F_q(\la)^\times\subset \G$ is isomorphic to $\F_{q^2}^\times$.
Hence $\la\mapsto \F_q(\la)^\times$ defines a map $\chi:P\to G$. As
before, it is easy to check that $\chi$ is well-defined and
surjective. We will show that $\chi$ is 2-to-1, which implies the
formula for $V_1$. Obviously $\la\neq \la'$ and these elements
generate the same subgroup in $\G$, as the canonical involution on
$D$ restricted to $F(\la)$ is equal to the Galois conjugation on
$F(\la)/F$. Since $\la$ and $\la'$ are the only elements in
$\F_q(\la)$ with the given characteristic polynomial, it is enough
to show that $\la$ and $\la'$ are not $\G$-conjugate. Suppose there
is $\gamma\in \G$ such that $\la'=\gamma \la \gamma^{-1}$. One
easily checks that $1, \la, \gamma, \gamma\la$ are linearly
independent over $F$, hence generate $D$. If $q$ is odd, then
$\la'=-\la$. If $q$ is even, then $\la'=\la+1$. Using this, one
easily checks that $\gamma^2$ commutes with $\la$, e.g., for $q$
odd:
$$
\gamma^2 \la \gamma^{-2}=\gamma \la'\gamma^{-1}=-\gamma \la
\gamma^{-1}=-\la'=\la.
$$
Hence $\gamma^2$ lies in the center of $D$, and therefore,
$\gamma^2=b\in \F_q^\times$. Looking at the relations between $\la$
and $\gamma$, we see that $D$ is isomorphic to $H(\xi, b)$. This
last quaternion algebra is split according to Corollary
\ref{cor-split2}, which leads to a contradiction.

The number of edges $E=\#\Ed(\G\bs\cT)$, ignoring the orientation,
is equal to
$$
E=(V_1+(q+1)V_{q+1})/2.
$$
By Euler's formula, $E+1=g(R)+V_1+V_{q+1}$. This implies the
expression for $V_{q+1}$, and finishes the proof of the theorem.
\end{proof}

The next theorem is the group-theoretic incarnation of Theorem
\ref{PropTV}.

\begin{thm}\label{thmGTI} Let $\G_\tor$ be the normal subgroup of $\G$ generated
by torsion elements.
\begin{enumerate}
\item $\G/\G_\tor$ is a free group on $g(R)$ generators.
\item If $\wp(R)=0$, then $\G_\tor=\F_q^\times$.
\item If $\wp(R)=1$, then the maximal finite order subgroups of $\G$
are isomorphic to $\F_{q^2}^\times$, and, up to conjugation, $\G$
has $2^{\# R-1}$ such subgroups.
\end{enumerate}
\end{thm}
\begin{proof}
By \cite[Cor. 1, p.55]{SerreT}, $\G/\G_\tor$ is the fundamental
group of the graph $\G\bs\cT$. The topological fundamental group of
any finite graph is a free group. Hence $\G/\G_\tor$ is a free
group. The number of generators of this group is equal to $\dim_\Q
H_1(\G\bs\cT, \Q)$.  This proves (1). Parts (2) and (3) follow from
the proofs of Proposition \ref{Prop2.1} and Theorem \ref{PropTV}.
\end{proof}

\begin{thm}\label{thm-tree}
$\G=\G_\tor$ if and only if one of the following holds:
\begin{enumerate}
\item $R=\{x,y\}$ and $\deg(x)=\deg(y)=1$. In this case, $\G$ has a
presentation
$$
\G\cong \langle \gamma_1,\gamma_2\ |\
\gamma_1^{q^2-1}=\gamma_2^{q^2-1}=1,\
\gamma_1^{q+1}=\gamma_2^{q+1}\rangle.
$$
\item $q=4$ and $R$ consists of the
four degree-$1$ places in $|C|-\infty$. In this case, $\G$ has a
presentation
$$
\G\cong \langle \gamma_1,\dots,\gamma_8\ |\ \gamma_1^{q^2-1}=\cdots
=\gamma_8^{q^2-1}=1,\ \gamma_1^{q+1}=\cdots=\gamma_8^{q+1}\rangle.
$$
\end{enumerate}
\end{thm}
\begin{proof} By Theorem \ref{thmGTI},
$\G=\G_\tor$ if and only if $g(R)=0$. From the formula for $g(R)$
one easily concludes that $g(R)=0$ exactly in the two cases listed
in the theorem. In Case (1), according to Theorem \ref{PropTV},
$V_1=2$ and $V_{q+1}=0$, so $\G\bs \cT$ is an edge:
\begin{center}
\begin{picture}(30,10)
\put(5,5){\circle*{2}}\put(5,5){\line(1,0){20}}
\put(25,5){\circle*{2}}
\end{picture}
\end{center}
In Case (2), $V_1=8$ and $V_{q+1}=2$, so $\G\bs \cT$ is the tree:
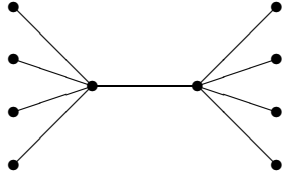
\begin{center}
\begin{picture}(60,40)
\put(5,5){\circle*{2}}\put(5,5){\line(1,1){15}}
\put(5,15){\circle*{2}}\put(5,15){\line(3,1){15}}
\put(5,25){\circle*{2}}\put(5,25){\line(3,-1){15}}
\put(5,35){\circle*{2}}\put(5,35){\line(1,-1){15}}
\put(20,20){\circle*{2}}\put(20,20){\line(1,0){20}}
\put(40,20){\circle*{2}}
\put(55,5){\circle*{2}}\put(55,5){\line(-1,1){15}}
\put(55,15){\circle*{2}}\put(55,15){\line(-3,1){15}}
\put(55,25){\circle*{2}}\put(55,25){\line(-3,-1){15}}
\put(55,35){\circle*{2}}\put(55,35){\line(-1,-1){15}}
\end{picture}
\end{center}
The knowledge of the quotient tree $\G\bs \cT$ allows to reconstruct
$\G$ \cite[I.4.4]{SerreT}: $\G$ is the graph of groups $\G\bs \cT$,
where each terminal vertex of $\G\bs \cT$ is labeled by
$\F_{q^2}^\times$, each non-terminal vertex is labeled by
$\F_q^\times$, and each edge is labeled by $\F_q^\times$ (the
monomorphisms $\G_y\to \G_{t(y)}$ are the natural inclusions
$\F_q^\times\hookrightarrow \F_{q^2}^\times$); see \cite[p.
37]{SerreT}. In other words, $\G$ is the amalgam of the groups
labeling the vertices of $\G\bs\cT$ along the subgroups labeling the
edges. The presentation for $\G$ follows from the definition of
amalgam; see \cite[I.1]{SerreT}.
\end{proof}

Theorem \ref{PropTV} allows to determine $\G\bs \cT$ in some other
cases, besides the case when $\G\bs \cT$ is a tree treated in
Theorem \ref{thm-tree}:

\begin{thm}\label{thm-hyp}
Suppose $R=\{x, y\}$ and $\{\deg(x), \deg(y)\}=\{1,2\}$. Then $\G\bs
\cT$ is the graph which has $2$ vertices and $q+1$ edges connecting
them:
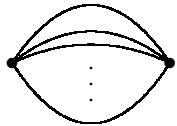
\begin{center}
\begin{picture}(40,25)
\qbezier(5,13)(20,35)(35,13)\qbezier(5,13)(20,25)(35,13)
\qbezier(5,13)(20,20)(35,13) \qbezier(5,13)(20,-10)(35,13)
\put(5,13){\circle*{2}}\put(35,13){\circle*{2}}
\put(20,12){\circle*{.7}}\put(20,9){\circle*{.7}}\put(20,6){\circle*{.7}}
\end{picture}
\end{center}
\end{thm}
\begin{proof}
From Theorem \ref{PropTV}, $V_1=0$ and $V_{q+1}=2$. This implies the
claim.
\end{proof}

The example in Theorem \ref{thm-hyp} is significant for arithmetic
reasons. Assume $q$ is odd. As is shown in \cite{PapHD}, the curve
$X^R$ is hyperelliptic if and only if $R=\{x, y\}$ and $\{\deg(x),
\deg(y)\}=\{1,2\}$. Thus, Theorems \ref{thm-hyp} and \ref{thmKM}
imply that the closed fibre over $\infty$ of the minimal regular
model over $\Spec(\cO)$ of a hyperelliptic $X^R$ consists of two
projective lines $\P^1_{\F_q}$ intersecting transversally at their
$q+1$ $\F_q$-rational points. This can be used to determine the
group of connected components $\Phi_\infty$ of the closed fibre of
the N\'eron model of the Jacobian of $X^R$ over $\Spec(\cO)$. Using
a result of Raynaud (cf. \cite[p. 283]{NM}), one obtains
$\Phi_\infty\cong \Z/(q+1)\Z$.


\section{Modular curves}\label{SecMC} $\GL_2(K)$ acts on Drinfeld's upper half-plane $\Omega$
by linear fractional transformations. As we discussed in the
introduction, the quotient $\G\bs \Omega$ is the underlying
rigid-analytic space of a smooth, projective curve $X^R$ over $K$.
The genus of this curve is computed in \cite{PapGenus} using
arithmetic methods, and it turns out to be equal to $g(R)$. The
theory of Mumford curves allows to construct a model of $X^R$ over
$\Spec(\cO)$.

\begin{thm}\label{thmKM}
There is a scheme $\cX^R$ over $\Spec(\cO)$ which is proper, flat
and regular, and whose generic fibre
$\cX^R_K:=\cX^R\times_\cO\Spec(K)$ is isomorphic to $X^R$. The
geometric special fibre
$\cX^R_{\bar{k}}:=\cX^R\times_\cO\Spec(\bar{k})$ is reduced,
connected and one dimensional and has at most ordinary double
points, where $\bar{k}$ denotes the algebraic closure of $k$. The
normalization of components of $\cX^R_{\bar{k}}$ are $k$-rational
curves, and the double points of $\cX^R_{\bar{k}}$ are $k$-rational
with two $k$-rational branches. Moreover, $\G\bs \cT$ is the dual
graph of $\cX^R_{\bar{k}}$: this means that the irreducible
components $E$ of $\cX^R_{\bar{k}}$ and the double points $x$ of
$\cX^R_{\bar{k}}$ are naturally in bijection with the vertices $v$
and the edges $\{y, \bar{y}\}$ of $\G\bs \cT$ respectively, such
that $x$ is contained in $E$ if and only if $v=o(y)$ or $t(y)$.
\end{thm}
\begin{proof}
This follows from \cite[Prop. 3.2]{Kurihara} after making the
following two observations: (1) since $\G$ acts without inversion of
$\cT$, the graph $(\G\bs\cT)^\ast$ in \cite{Kurihara} is
$(\G\bs\cT)$ itself; (2) for any $y\in \Ed(\G\bs \cT)$ the image of
the stabilizer $\G_y$ is trivial in $\PGL_2(\Fi)$ by Proposition
\ref{Prop2.1}, so the length of $y$ (in Kurihara's terminology) is
$1$.
\end{proof}

\begin{rem}
$\cX^R$ is not necessarily the minimal regular model of $X^R$ over
$\Spec(\cO)$; one obtains the minimal model by removing the terminal
vertices from $\G\bs\cT$.
\end{rem}

Next, we recall Jordan-Livn\'e's geometric version of Hensel's
lemma:

\begin{thm}\label{thmJL}
Let $\cX$ be a proper, flat and regular scheme over $\cO$. Let
$\cX_K:=\cX\times_\cO\Spec(K)$ be the generic fibre and $\cX_k:=\cX
\times_\cO\Spec(k)$ be the special fibre. Then $\cX_K$ has a
$K$-rational point if and only if $\cX_k$ has a smooth $k$-rational
point.
\end{thm}
\begin{proof}
See \cite[Lem. 1.1]{JL}.
\end{proof}

\begin{thm}\label{thmFipoints}
$X^R(K)\neq \emptyset$ if and only if $\wp(R)=1$.
\end{thm}
\begin{proof}
By Theorem \ref{thmKM}, $X^R$ is the generic fibre of $\cX^R$ which
is proper, flat and regular over $\Spec(\cO)$. Hence by Theorem
\ref{thmJL}, $X^R$ has a $K$-rational point if and only if $\cX^R_k$
has a smooth $k$-rational point. The irreducible components of
$\cX^R_k$ are $\P^1_k$'s intersecting transversally at $k$-rational
points, and exactly two components pass through a singularity (note
that the irreducible components do not have self-intersections, as
$\G\bs\cT$ has no loops.) Since $\#\P^1_k(k)=q+1$, $\cX^R_k$ has a
smooth $k$-rational point if and only if there is an irreducible
component whose corresponding vertex in $\G\bs\cT$ has degree less
than $q+1$. By Theorem \ref{PropTV}, $\G\bs\cT$ has a vertex of
degree less than $q+1$ if and only if $\wp(R)=1$.
\end{proof}


\section{Explicit generators}\label{SecExplU} Over $\Q$ one knows not only the cases
when $\G^d$ is generated by torsion elements, but also the explicit
description of the generators of $\G^d$ in terms of a basis of $D$;
cf. \cite[p. 92]{AB} or \cite{KV}. For example, $\G^6$ is isomorphic
to the subgroup of $\SL_2(\R)$ generated by
\begin{align*}
\gamma_1=\frac{1}{2}\begin{pmatrix} \sqrt{2} & 2-\sqrt{2} \\
-6-3\sqrt{2} & -\sqrt{2}\end{pmatrix}&,\quad
\gamma_2=\frac{1}{2}\begin{pmatrix} \sqrt{2} & -2+\sqrt{2} \\
6+3\sqrt{2} & -\sqrt{2}\end{pmatrix}\\
\gamma_3=\frac{1}{2}\begin{pmatrix} 1 & 1 \\
-3 & 1\end{pmatrix}&,\quad
\gamma_4=\frac{1}{2}\begin{pmatrix} 1 & 3-2\sqrt{2} \\
-9-6\sqrt{2} & 1\end{pmatrix}
\end{align*}
which have orders $4$, $4$, $6$, $6$, respectively.

In this section we will find explicit generators for $\G$ in Case
(1) of Theorem \ref{thm-tree}. As a consequence of our calculations,
we will also obtain in this case a direct proof of Theorem
\ref{PropTV} for odd $q$.

First, we explicitly describe $D$ in terms of generators and
relations when $\wp(R)=1$, and then describe a maximal $A$-order in
$D$. For each $x\in |C|-\infty$, denote by $\fp_x\lhd A$ the
corresponding prime ideal of $A$. Let $\fr$ be the monic generator
of the ideal $\prod_{x\in R}\fp_x$ (this is the
\textit{discriminant} of $D$).

\begin{lem}
Assume $\wp(R)=1$. If $q$ is odd, let $\xi\in \F_{q}$ be a fixed
non-square. If $q$ is even, let $\xi\in \F_{q}$ be a fixed element
such that $\Tr_{\F_q/\F_2}(\xi)\neq 0$. Then $H(\xi,\fr)\cong D$.
\end{lem}
\begin{proof}
It is enough to check that $H$ is ramified exactly at the places in
$R$. For this one can use the same argument as in the proof of
Corollary \ref{cor-split2}. Reducing the quadratic form $Q$ mod
$\fp_v$ for every $v\in |C|-\infty$, one easily checks that $H(\xi,
\fr)$ is ramified at every place in $R$ and is split at every place
in $|C|-R-\infty$: Note that $x^2-\xi y^2=0$ (resp. $x^2+xy+\xi
y^2=0$) has no non-trivial solutions over $\F_v$, $v\in R$, when $q$
is odd (resp. even), since the corresponding quadratic has no roots
over $\F_q$ due to the choice of $\xi$ and $[\F_v:\F_q]$ is odd by
assumption. Finally, $H$ is automatically split at $\infty$ since
the number of places where a quaternion algebra ramifies must be
even.
\end{proof}

From now on we assume that $\wp(R)=1$ and identify $D$ with $H(\xi,
\fr)$.

\begin{lem}\label{lem3.5}
The free $A$-module $\La$ in $D$ generated by
$$
x_1=1,\quad x_2=i,\quad x_3=j,\quad x_4=ij
$$
is a maximal order.
\end{lem}
\begin{proof}
It is obvious that $\La$ is an order. To show that it is maximal we
compute its discriminant, i.e., the ideal of $A$ generated by
$\det(\Tr(x_ix_j))_{ij}$. When $q$ is odd
$$
\det(\Tr(x_ix_j))_{ij}=\det \begin{pmatrix} 2 & 0 & 0 & 0\\
0 & 2\xi & 0 & 0\\
0 & 0 & 2\fr & 0\\
0 & 0 & 0 & -2\xi \fr
\end{pmatrix}=-16\xi^2\fr^2,
$$
when $q$ is even
$$
\det(\Tr(x_ix_j))_{ij}=\det \begin{pmatrix}
0 & 1 & 0 & 0\\
1 & 0 & 0 & 0\\
0 & 0 & 0 & \fr\\
0 & 0 & \fr & 0
\end{pmatrix}=\fr^2.
$$
In both cases, the discriminant of $\La$ is $\prod_{x\in R}\fp_x^2$,
and this implies that $\La$ is maximal, cf. \cite[pp.
84-85]{Vigneras}.
\end{proof}

From Lemma \ref{lem3.5} we see that finding the elements
$\la=a+bi+cj+dij\in D$ which lie in $\G$ is equivalent to finding
$a,b,c,d,\in A$ such that
$$
(a^2-\xi b^2)-\fr(c^2-\xi d^2) \in \F_q^\times \quad \text{ if $q$
is odd}
$$
$$
a^2+ab+\xi b^2+\fr(c^2+cd+\xi d^2)\in \F_q^\times \quad \text{ if
$q$ is even},
$$
(this is the condition $\Nr(\la)\in \F_q^\times$ written out
explicitly). We are particularly interested in torsion elements of
$\G$, hence instead of looking for general units, we will try to
find elements in $\La$ which are algebraic over $\F_q$ (such
elements automatically lie in $\La^\times$). Since the actual
calculations differ for ``$q$ even'' and ``$q$ odd'', we have to
treat these cases separately, but first we make the following
simplifying observation. If $R$ consists of two rational places,
then without loss of generality we can assume $\fr=T(T-1)$. Indeed,
$F$ is the function field of $\P^1_{\F_q}$, so $\Aut(F/\F_q)\cong
\PGL_2(\F_q)$, where the matrix $\begin{pmatrix} x & y\\ z &
w\end{pmatrix}$ acts by $T\mapsto \frac{xT+y}{zT+w}$. It is
well-known that there is a linear fractional transformation which
fixes $\infty$ and maps any two given $\F_q$-rational points of
$\P^1_{\F_q}$ to $0$ and $1$. The action of $\Aut(F/\F_q)$ does not
affect the structure of $\G\bs \cT$ (but of course, quaternion
algebras ramified at different sets of places are not isomorphic as
$F$-algebras).

\subsection*{q is odd} Consider the equation $\gamma^2=\xi$ in $\La$.
If we write $\gamma=a+bi+cj+dij$, then
$$
\gamma^2= a^2+b^2\xi +c^2\fr-d^2\xi\fr+2(abi+acj+adij).
$$
Therefore, $\gamma^2=\xi$ is equivalent to
\begin{equation}\label{eq-tuodd}
a=0\quad \text{and} \quad b^2\xi  +c^2\fr-d^2\xi\fr=\xi.
\end{equation}
A possible solution is $b=1$, $d=c=0$. This gives the obvious
$\theta_1=i$ as a torsion unit.

Now let $\fr=T(T-1)$. One easily checks that
$$
b=2T-1,\quad c=0, \quad d=2
$$
satisfies (\ref{eq-tuodd}). Thus, $\theta_2=(2T-1)i+2ij$ is a
torsion unit.

Next, we study the action of $\theta_1, \theta_2$ on $\cT$. Since
$\fr=T(T-1)$ has even degree and is monic, $\sqrt{\fr}\in K$. The
map
$$
i\mapsto \begin{pmatrix} 0 & 1\\ \xi & 0 \end{pmatrix}, \quad
j\mapsto
\begin{pmatrix} \sqrt{\fr} & 0\\ 0 & -\sqrt{\fr}\end{pmatrix},
$$
defines an embedding of $D$ into $\M_2(K)$.

Consider $L:=F(\sqrt{\fr})$. The place $\infty$ splits in $L$.
Denote by $\infty_1$ and $\infty_2$ the two places of $L$ over
$\infty$. The integral closure of $A$ in $L$ is
$\cA:=A[\sqrt{\fr}]$. Let $\pi:=((2T-1)+2\sqrt{\fr})\in \cA$. Since
$\pi^{-1}=(2T-1)-2\sqrt{\fr}$, $\pi$ is a unit in $\cA$. Since $\pi$
is not a constant, its valuations are non-zero at $\infty_1$ and
$\infty_2$. If $n:=\ord_{\infty_1}(\pi)$, then
$$
\min(n, -n)=\ord_{\infty_1}(\pi+\pi^{-1})=\ord_{\infty_1}(4T-2)=-1.
$$
This implies that we can extend the valuation at $\infty$ to $L$ so
that $\pi$ is a uniformizer of $F_\infty$. We conclude that
$\theta_2$ acts on $\cT$ as the matrix
$$
\theta_2=i((2T-1)+2j)=\begin{pmatrix} 0 & 1\\ \xi & 0
\end{pmatrix} \begin{pmatrix} \pi  & 0 \\ 0 & \pi^{-1}
\end{pmatrix}=\begin{pmatrix} 0 & \pi^{-1}\\ \xi\pi & 0
\end{pmatrix}.
$$

Let $v, w\in \Ver(\cT)$ be $v=[e_1, e_2]$ and $w=[\pi e_1, e_2]$.
Now
$$
\theta_1 \cdot v = [e_2, \xi e_1] = [\xi e_1, e_2]=[e_1, e_2]=v
$$
and
$$
\theta_2\cdot w=[e_2, \xi\pi e_1]=[\xi\pi e_1, e_2]=[\pi e_1,
e_2]=w.
$$
Thus, we found two adjacent vertices in $\cT$ and elements in their
stabilizers which are not in $\F_q^\times$. Proposition
\ref{Prop2.1} implies that $\G_v\cong\F_{q^2}^\times$ and
$\G_w\cong\F_{q^2}^\times$. By Corollary \ref{cor2.2}, the images of
$v$ and $w$ in $\G\bs \cT$ are adjacent terminal vertices. This
implies that $\G\bs\cT$ is an edge, and proves Theorem \ref{PropTV}
in the case when $R=\{x, y\}$ and $\deg(x)=\deg(y)=1$. By fixing
generators $\gamma_1, \gamma_2$ of the finite cyclic groups
$\F_q(\theta_1)^\times=\langle\gamma_1\rangle$ and
$\F_q(\theta_2)^\times=\langle\gamma_2\rangle$, one obtains two
torsion elements which generate $\G$.

\begin{example}
Let $q=3$. Then $\xi=-1$ and $1-\theta_i$ generates
$\F_q(\theta_i)^\times$. Hence $\G$ is isomorphic to the subgroup of
$\GL_2(K)$ generated by the matrices
$$
\gamma_1=\begin{pmatrix} 1 & 2 \\ 1 & 1\end{pmatrix}\quad
\text{and} \quad \gamma_2=\begin{pmatrix} 1 & (T+1)+2\sqrt{T(T-1)} \\
2(T+1)+2\sqrt{T(T-1)} & 1\end{pmatrix}
$$
both of which have order $8$ and satisfy $\gamma_1^4=\gamma_2^4=-1$.
\end{example}

\subsection*{q is even} Here we will look for solutions of the
equation $\gamma^2+\gamma + \xi=0$ in $\La$. Again writing
$\gamma=a+bi+cj+dij$,
\begin{align*}
\gamma^2+\gamma +\xi =& (a+a^2+b^2\xi+c^2\fr+cd\fr+d^2\xi\fr+\xi)\\
&+b(b+1)i+c(b+1)j+d(b+1)ij.
\end{align*}
Hence $\gamma^2+\gamma +\xi =0$ if and only if
\begin{equation}\label{eq-tu}
b=1 \quad \text{and} \quad a+a^2+\fr(c^2+cd+d^2\xi)=0.
\end{equation}
An obvious solution is $a=c=d=0$, which gives $\theta_1=i$ as a
torsion unit.

Now let $\fr=T(T+1)$. Then
$$
a=T, \quad c=1,\quad d=0
$$ satisfy (\ref{eq-tu}), so $\theta_2=T+i+j$
is a torsion unit.

When $q$ is even, a technical complication arises in the study of
the action of $\G$ on $\cT$. This is due to the fact that neither
$F(i)$ nor $F(j)$ embed into $K$ (in the extension $F(i)/F$ the
place $\infty$ remains inert, and $\infty$ ramifies in the extension
$F(j)/F$, since $j$ is non-separable over $F$.) In particular, our
choice of generators of $D$ does not provide an easy explicit
embedding of $\G$ into $\GL_2(K)$ (although, we can embed $\G$ into
$\GL_2(\F_{q^2}K)$ via
$$
i\mapsto \begin{pmatrix} \kappa & 0\\0 & \kappa+1\end{pmatrix},
\quad j\mapsto
\begin{pmatrix} 0 & 1\\ \fr & 0\end{pmatrix},
$$
where $\kappa\in \F_{q^2}$ is a solution of $x^2+x+\xi=0$). Instead,
we appeal to Theorem \ref{thmGTI}. A tedious calculation shows that
$\F_q(\theta_1)^\times$ and $\F_q(\theta_2)^\times$ are not
conjugate in $\G$ (in fact, it is enough to check that $\theta_1$ is
not a $\G$-conjugate of $\theta_2$ or $\theta_2+1$). Hence by
Theorem \ref{thmGTI} these subgroups generate $\G$. Again, choosing
generators $\gamma_1$ and $\gamma_2$ of these cyclic groups, one
obtains two torsion elements which generate $\G$.

\vspace{0.1in}

We finish with a remark on Case (2) of Theorem \ref{thm-tree}: $q=4$
and $\fr=T^4+T$. We can write down some solutions of (\ref{eq-tu})
in this case. Let $c,d\in \F_4$ be not both zero. Then
$$
\alpha:=c^2+cd+\xi d^2\in \F_4^\times.
$$
Let $s\in \F_4^\times$ be such that $s^2=\alpha$ (always exists and
is unique). Let $a=sT^2+s^2T+m$, where $m=0,1$. Then, since $s^4=s$
and $m^2=m$,
$$
a^2+a=(\alpha T^4+sT^2+m)+(sT^2+\alpha T+m)=\alpha(T^4+T).
$$
Now it is clear that $\theta=a+i+cj+dij$, with $a,c,d$ as above,
satisfies (\ref{eq-tu}). Nevertheless, it seems rather challenging
to find explicitly enough torsion units which will generate the $8$
non-conjugate subgroups of $\G$ isomorphic to $\F_{16}^\times$.


\subsection*{Acknowledgments}  I thank E.-U. Gekeler and A. Schweizer  for
useful discussions. The article was written while I was visiting the
Department of Mathematics of Saarland University.



\end{document}